\documentclass[acmsmall]{acmart}

\usepackage{listings}
\usepackage{algorithm}
\usepackage{algpseudocode}
\usepackage{mathdots}
\AtBeginDocument{%
  \providecommand\BibTeX{{%
    \normalfont B\kern-0.5em{\scshape i\kern-0.25em b}\kern-0.8em\TeX}}}

\setcopyright{acmcopyright}
\copyrightyear{2024}
\acmYear{2024}
\acmDOI{10.1234/1122334.1122345}

\acmJournal{JACM}
\acmVolume{0}
\acmNumber{0}
\acmArticle{0}
\acmMonth{0}

\usepackage{cleveref}
\usepackage{todonotes}

\usepackage{physics}
\usepackage{bbold}
\usepackage{siunitx}
\newtheorem{remark}{Remark}[section]

\newcommand{\pythOS}{{\tt pythOS}}
\newcommand{\Irksome}{{\tt Irksome}}
\newcommand{\Firedrake}{{\tt Firedrake}}

\newcommand{\Dt}{\Delta t}
\newcommand{\Dtn}{\Dt_n}
\newcommand{\Dx}{\Delta x}
\newcommand{\Dy}{\Delta y}
\newcommand{\tf}{t_{\text{f}}}

\newcommand{\FF}{\mathbf{F}}
\newcommand{\F}{\boldsymbol{\mathcal{F}}}
\newcommand{\Fl}[1]{\boldsymbol{\mathcal{F}}^{[#1]}}
\newcommand{\Gl}[1]{\boldsymbol{\mathcal{G}}^{[#1]}}
\newcommand{\YY}{\mathbf{Y}}
\newcommand{\yy}{\mathbf{y}}

\newcommand{\ssl}{s^{[\ell]}}
\newcommand{\sslp}{s^{[\ell']}}
\newcommand{\cc}[1]{\mathbf{c}^{[#1]}}
\newcommand{\bb}[1]{\mathbf{b}^{[#1]}}
\newcommand{\AAA}[1]{\mathbf{A}^{[#1]}}
\newcommand{\ty}{\tilde{\yy}}

\newcommand{\aaalpha}[2]{\alpha_{#1}^{[#2]}}
\newcommand{\mrislowstage}{s^{[S]}}
\newcommand{\Nop}{N}
\newcommand{\Zero}{\boldsymbol{0}}
\newcommand{\GGamma}[1]{\mathbf{\Gamma}^{\{#1\}}}
\newcommand{\OOmega}[1]{\mathbf{\Omega}^{\{#1\}}}

\newcommand{\CLTTwo}{\text{Complex\_Lie\_Trotter\_2}}
\newcommand{\CLTThree}{\text{Complex\_Lie\_Trotter\_3}}
\newcommand{\complexde}{\mathcal{C}}
\newcommand{\realde}{\mathcal{R}}

\newcommand{\numpy}{\texttt{numpy}}
\newcommand{\scipy}{\texttt{scipy}}

\newcommand{\atol}{\texttt{atol}}
\newcommand{\rtol}{\texttt{rtol}}

\newcommand{\one}{\mathbb{1}}
\newcommand{\Ny}{\tilde{N}}
\newcommand{\faststepind}{m}

\begin{document}

\title{\pythOS: A Python library for solving IVPs by operator splitting}

\author{Victoria Guenter}
\authornotemark[1]
\email{v.guenter@usask.ca}
\affiliation{%
  \institution{University of Saskatchewan}
  \streetaddress{Department of Computer Science, 110 Science Place}
  \city{Saskatoon}
  \state{Saskatchewan}
  \country{Canada}
  \postcode{S7N 5C9}
}
\author{Siqi Wei}
\authornotemark[2]
\email{siqi.wei@usask.ca}
\affiliation{%
  \institution{University of Saskatchewan}
  \streetaddress{Department of Mathematics and Statistics, 110 Science Place}
  \city{Saskatoon}
  \state{Saskatchewan}
  \country{Canada}
  \postcode{S7N 5C9}
}
\author{Raymond J.~Spiteri}
\authornotemark[1]
\email{spiteri@cs.usask.ca}
\orcid{0000-0002-3513-6237}
\affiliation{%
  \institution{University of Saskatchewan}
  \streetaddress{Department of Computer Science, 110 Science Place}
  \city{Saskatoon}
  \state{Saskatchewan}
  \country{Canada}
  \postcode{S7N 5C9}
}

\renewcommand{\shortauthors}{Guenter, Wei, and Spiteri}

\begin{abstract}
  Operator-splitting methods are widespread in the numerical solution
  of differential equations, especially the initial-value problems in
  ordinary differential equations that arise from a method-of-lines
  discretization of partial differential equations. Such problems can
  often be solved more effectively by treating the various terms
  individually with specialized methods rather than simultaneously in
  a monolithic fashion. This paper describes \pythOS, a Python
  software library for the systematic solution of differential
  equations by operator-splitting methods. The functionality of
  \pythOS\ focuses on fractional-step methods, including those with
  real and complex coefficients, but it also implements additive
  Runge--Kutta methods, generalized additive Runge--Kutta methods, and
  multi-rate, and multi-rate infinitesimal methods. Experimentation
  with the solution of practical problems is facilitated through an
  interface to the \Firedrake\ library for the finite element spatial
  discretization of partial differential equations and further
  enhanced by the convenient implementation of exponential
  time-integration methods and fully implicit Runge--Kutta methods
  available from the \Irksome\ software library.  The functionality of
  \pythOS\ as well as some less generally appreciated aspects of
  operator-splitting methods are demonstrated by means of examples.
\end{abstract}


\maketitle

\section{Introduction}

Differential equations play a fundamental role in many fields of
science and engineering. Due to the size or complexity of specific
systems being modelled, it may not be optimal to solve the
differential equations monolithically, i.e., using one numerical
method to solve the entire system, as is standard practice. For
example, different parts of the differential equation may benefit from
specialized solvers. In such cases, operator-splitting methods are
often used as a divide-and-conquer approach to make numerical
simulations feasible or more efficient.

The focus of this paper is on the numerical solution of initial-value
problems (IVPs) in ordinary differential equations (ODEs).  It is
often natural for the right-hand side of an ODE to be comprised of
multiple contributing terms, taking the $\Nop$-additive form
\begin{equation}
	\label{eq:N-ODE}
	\dv{\yy}{t} = \F(t,\yy) = \sum_{\ell=1}^{\Nop} \Fl{\ell} (t, \yy),
    \qquad \yy(0) = \yy_0.
\end{equation}
This form is common when applying the method-of-lines discretization
to partial differential equations, especially those that model
multi-physics phenomena. For example, in models of combustive flows,
there are contributions from advection, reaction, and
diffusion. This form also includes the \textit{partitioned} form,
\begin{align*}
  \dv{\yy_1}{t} &= \FF_{1}(\yy_1,\yy_2,\dots,\yy_{\Ny}), \\
  \dv{\yy_2}{t} &= \FF_{2}(\yy_1,\yy_2,\dots,\yy_{\Ny}), \\
                & \ \vdots \\
    \dv{\yy_{\Ny}}{t} &= \FF_{\Nop}(\yy_1,\yy_2,\dots,\yy_{\Ny}),
\end{align*}
by defining
$\Fl{\ell}(\yy) =
(\Zero,\dots,\Zero,\FF_\ell((\yy_1,\yy_2,\dots,\yy_{\Ny})),
\Zero,\dots,\Zero )$, $\ell=1,2,\dots,\Nop$. Although, in principle,
the number of $\yy$ components, $\Ny$, does not equal to the number of
partitions, $\Nop$, in practice, $\Ny$ is often equal to $\Nop$.

There are a variety of techniques that may be applied to leverage the
additive structure of~\cref{eq:N-ODE}.  One of the most common of
these is the class of \emph{fractional-step methods},
e.g.,~\cite{Yanenko1971}. Fractional-step methods can generally be
thought of as operator-splitting methods where the operators
$\Fl{\ell}(t,\yy)$ (or \emph{fractions} of the right-hand side) are
integrated separately (potentially with different integrators) over
sub-steps of the time-step $[t_{n},t_{n+1}]$, and the solution at the
end of the time-step is formed from a combination of these
sub-integrations.  Such methods are particularly popular for
simulations of multi-physics problems, where individual
sub-integrations may be performed in a black-box fashion using
separate libraries --- so-called
\textit{co-simulations}~\cite{Gomes2018-coSimulation}.

Another highly common operator-splitting method is the class of
implicit-explicit (IMEX) methods~\cite{Ascher1995, ars1997}, which are
2-additive methods that treat one of the right-hand side terms (e.g.,
the stiff term) with an implicit method and the other (e.g., the
non-stiff term) with an explicit method. IMEX methods have often met
with success because the marginal loss in stability compared to a
fully implicit method is more than compensated for by the reduction in
computational cost per step.  It can be shown that fractional-step
Runge--Kutta methods are a special case of additive Runge--Kutta
methods; i.e., the additive Butcher tableau representation has a particular
structure~\cite{spiteri_wei_FSRK} that reflects a specific form of coupling
between operators, namely that the influence of one operator on
another only occurs at the initial condition of a sub-integration.

Other operator-splitting methods include $N$-additive Runge--Kutta
(ARK) methods \cite{kennedy2003}, generalized additive Runge--Kutta
(GARK) methods \cite{sandu2015}, multi-rate (MR) methods
\cite{knoth1998, wensch2009}, and multi-rate infinitesimal (MRI)
methods~\cite{wensch2009,Sandu2019,Chinomona2021}.  Here, we use the
term ``operator-splitting methods'' to describe any method that
involves some kind of non-monolithic time discretization. However, some
authors use this general term to specifically describe what we call
``fractional-step methods'' in this paper, e.g.,~\cite{macnamara2016}.

Despite their ubiquity in practice, most implementations of
operators-splitting methods in multi-physics applications are made on
an ad-hoc, one-off basis for a given problem. Such implementations are
often low order (Lie--Trotter or Strang
splitting)~\cite{Godunov1959,strang1963}. Reasons for this may include
difficulty to determine which splitting method to use and possible
backward-in-time integration required for high-order fractional-step
methods \cite{sheng1989,suzuki1991,goldman1996}. Also, many such
applications are based on discretizations of partial differential
equations (PDEs). Such problems require discretization in both space
and time, and many finite element libraries do not offer systematic
support for users to implement their own time stepping.

There are software packages such as ARKODE \cite{arkode}, PESTc/TS
\cite{Abhyankar2018pestc}, and Trilinos/Tempus \cite{OberTrilinos}
that have been developed to flexibly perform time integration using
additive Runge--Kutta methods.  The ARKODE library supports
Runge--Kutta, IMEX-Runge--Kutta, and 2- or 3-additively split MRI
methods of various orders, including adaptive Runge--Kutta and IMEX
Runge--Kutta methods, as well as user defined methods. PESTc and
Trilinos both support 2-additive IMEX Runge--Kutta methods, and PESTc
supports explicit multi-rate Runge--Kutta methods.



%

The \pythOS\ library focuses particularly on facilitating the
exploration of the fractional-step methods.  The use of a
fractional-step method to solve a given problem generally involves
three main steps:
\begin{enumerate}
	
\item \label{step:split} Write the differential equation in the
  form~\cref{eq:N-ODE} with the right-hand side split into $\Nop$
  additive operators (or sub-systems).
	
\item \label{step:sub-integ} Integrate the sub-systems so formed with
  appropriate sub-integration methods.
	
\item \label{step:couple} Combine the solutions of the sub-systems to
  obtain a solution to the original problem.
\end{enumerate}

How to perform Step~\ref{step:split} in an effective manner is
generally an open question and highly problem-dependent. Some
literature \cite{Preuss2022,Wei2024} compared splitting strategies
such as splitting statically based on physical process or splitting
dynamically by dynamical linearization. Moreover, how many operators
should the right-hand side be split into is also a non-trivial
problem. Although traditional methods focus on 2-split methods,
researchers have started to look at $\Nop$-split methods for
$\Nop \geq 3$ \cite{spiteri2023_3split,MaraBehRaedAli2023}.  Applying
specialized sub-integration methods to each sub-system in
Step~\ref{step:sub-integ} is one of the main advantages of using
operator-splitting methods. For example, one can apply suitable
Runge--Kutta methods based on the stiffness of the operators, use the
exact solution of the sub-system if available, or use adaptive or
exponential methods for highly accurate solutions. The coupling
associated with Step~\ref{step:couple} typically has a significant
impact on the stability, accuracy, and computational cost of general
operator-splitting methods.

One goal of the \pythOS\ library is to allow researchers and users to
conveniently compare different operator-splitting methods and
experiment with various possibilities for each of the three steps:
splitting strategies, sub-integrations, and splitting schemes. The
\pythOS\ library allows the right-hand side function to be split into
any number of sub-systems and integrated in any order within a given
stage.  Integrations of sub-systems can be done with (adaptive)
Runge--Kutta methods, exponential methods, user-provided analytical
solutions, or adaptive multistep methods defined by external tools
(e.g., SUNDIALS and \scipy). In addition to existing methods, users
can easily implement their own fractional-step, ARK, GARK, MR, and MRI
methods.  Although the \pythOS\ library is designed to solve ODEs,
PDEs can be conveniently transformed to ODEs by the method of lines
with spatial discretization using \Firedrake.

The remainder of this paper is organized as follows.  In section
\ref{sec:background}, we present the mathematical framework for the
solvers implemented.  In section \ref{sec:structure}, we present the
structure of the \pythOS\ library.  In section \ref{sec:examples}, we
present some examples highlighting the various solvers and
non-standard operator-splitting options available. The coefficients of
all splitting methods used are provided in
appendix~\ref{sec:os_methods_coeff}. In section \ref{sec:conclusions},
we summarize our results and offer some conclusions.

\section{Mathematical Background}
\label{sec:background}

In this section, we introduce some of the mathematical background
behind the operator-splitting methods implemented in
\pythOS. Different splitting methods usually treat the operators
$\Fl{\ell}$ with specialized methods and combine the solutions through
a coupling mechanism. In \cref{subsec:FS}, we introduce the
fractional-step methods that in principle allow arbitrary
sub-integrators.  In \cref{subsec:ARK} and \cref{subsec:GARK}, we
introduce two families of splitting methods that use Runge--Kutta
sub-integrators. In \cref{subsec:MR}, we introduce the multi-rate and
multi-rate infinitesimal methods that specialize in treating
multi-physics problems with sub-systems that evolve on different
time scales.

\subsection{Fractional-Step methods}
\label{subsec:FS}
In this section, we introduce the fractional-step methods as defined
in \cite{hairer2006}.  After splitting the original ODE into $\Nop$
sub-systems as in~\cref{eq:N-ODE}, fractional-step methods integrate
one sub-system at a time. More explicitly, an $s$-stage, $\Nop$-split
fractional-step method is determined by the fractional-step
coefficients
$\{\aaalpha{k}{\ell} \}_{k=1,2,\dots, s}^{\ell = 1,2,\dots, \Nop}$. At
each stage $k$, let $\Phi_{\aaalpha{k}{\ell} \Dt}^{[\ell]}$ be the integrator applied
to sub-system $\ell$, 
\begin{equation}
\label{eq:sub-sys}
\dv{\yy^{[\ell]}}{t} = \Fl{\ell}(t,\yy^{[\ell]}),  
\end{equation}
over the step-size $\aaalpha{k}{\ell} \Dt$.
Let $\yy_n$ be the numerical solution of \cref{eq:N-ODE} at $t=t_n$. The
fractional-step method solves the additive system \cref{eq:N-ODE} over
the interval $[t_n, t_{n+1}]$ as outlined in
\cref{alg:general-os-step}.
\begin{algorithm}[!hbtp]
	\caption{\label{alg:general-os-step} Algorithm for taking a single
		step for a general fractional-step method.}
	
	\begin{algorithmic}[1]
		\Require $\aaalpha{k}{\ell} $, $\Fl{\ell}$, $\Phi_{\aaalpha{k}{\ell} \Dt_n}^{[\ell]}$, $t_n$, $t_{n+1}$, $\yy_n$.
		
		\State $\Dt_n=t_{n+1}-t_n$
		\State $t^{[\ell]}=t_n$ for $\ell=1,2,\ldots,\Nop$
		\State $\ty_0 = \yy_n$
		
		\For {$k=1$ to $s$}
		
		\For {$\ell=1$ to $\Nop$}
		\State $\alpha = \alpha_{k}^{[\ell]}$
		
		\State Solve
        $\displaystyle \dv{\ty^{[\ell]}}{t} = \Fl{\ell}\left(t,
          \ty^{[\ell]}\right) $, $\ty(t^{[\ell]})=\ty_0$, for
        $t \in \left[t^{[\ell]},t^{[\ell]}+\alpha\Dt_n\right]$ using
        $\Phi_{\alpha \Dt}^{[\ell]}$
		
		\State $t^{[\ell]} = t^{[\ell]} + \alpha \Dt_n$
		
		\State $\ty_0  =  \ty^{[\ell]}(t^{[\ell]})$
		
		\EndFor
		\EndFor
		\State \textbf{return} $\yy_{n+1} = \tilde\yy_0$
	\end{algorithmic}
\end{algorithm}


In the case where each sub-system is solved using one step of a
Runge--Kutta method, the overall fractional-step method is defined as
a fractional-step Runge--Kutta (FSRK) method
\cite{spiteri_wei_FSRK}. However, the sub-integration methods in
\pythOS\ are not limited to Runge--Kutta methods in this way. In
\pythOS, the sub-systems can be solved arbitrarily, using for example
adaptive Runge--Kutta methods, adaptive multi-step methods,
exponential methods, or exact solutions. It is also possible to use a
different sub-integration method for each operator and any stage.



\subsection{Additive Runge--Kutta methods}
\label{subsec:ARK}

In this section, we introduce ARK methods as described in
\cite{Cooper1980,kennedy2003}. After splitting an ODE into $\Nop$
sub-systems, each operator $\Fl{\ell}(t,\yy)$ is discretized using an
$s$-stage Runge--Kutta method. The resulting $s$-stage, $N$-additive
Runge--Kutta method advances the solution of an initial-value problem
\cref{eq:N-ODE} from $t_{n}$ to $t_{n+1}$ as follows:
\begin{align*}
	\mathbf{Y}_i & = \mathbf{y}_{n} +\Delta t \sum_{\ell=1}^{\Nop}\sum_{j=1}^{s} a_{ij}^{[\ell]}\Fl{\ell}(t_n+c_j^{[\ell]}\Delta t, \mathbf{Y}_j),\enskip i=1,2, \dots,s , \\ 
	\mathbf{y}_{n+1} & = \mathbf{y}_n +\Delta t\sum_{\ell=1}^{\Nop} \sum_{i=1}^{ s}b_i^{[\ell]}\Fl{\ell}(t_n+c_i^{[\ell]}\Delta t, \mathbf{Y}_i),\\
\end{align*}
where $a_{ij}^{[\ell]}$, $b_i^{[\ell]}$, and $c_j^{[\ell]}$ are the
coefficients of the Runge--Kutta method applied to operator
$\Fl{\ell}(t,\yy)$.  The Butcher tableau for ARK methods can be
written as \cite{sandu2015}
\begin{equation}\label{FSRKeq:arktab}
	\begin{array}{c|c|c|c|c|c|c|c}
		\mathbf{c}^{[1]}&\mathbf{c}^{[2]}&\cdots&\mathbf{c}^{[\Nop]}&\mathbf{A}^{[1]}&\mathbf{A}^{[2]}&\cdots &\mathbf{A}^{[\Nop]}\\
		\hline
		&&&&\mathbf{b}^{[1]}&\mathbf{b}^{[2]}&\cdots&\mathbf{b}^{[\Nop]}\\
	\end{array} ,
\end{equation}
where $\mathbf{A}^{[\ell]}$, $\mathbf{b}^{[\ell]}$,
$\mathbf{c}^{[\ell]}, \enskip \ell=1,2,\dots, \Nop$, are the Butcher
tableau quantities of the Runge--Kutta method associated with
operator $\ell$. The most common form of ARK is the 2-additive IMEX
method, where one of the Runge--Kutta methods is singly diagonally
implicit (and usually with good stability properties) and the other is
explicit.

\subsection{Generalized additive Runge--Kutta methods}
\label{subsec:GARK}

One generalization of the ARK methods allows each operator to be
solved by Runge--Kutta methods with different numbers of stages. 
This class of methods is known as the generalized additive
Runge--Kutta (GARK) methods as defined in \cite{sandu2015}.  
In particular, let $\ssl$ be the stage number of the Runge--Kutta method applied to operator $\Fl{\ell}$. Let
$a_{ij}^{[\ell',\ell]}$ and $b_j^{[\ell]}$,
$i=1, 2,\dots,\sslp,\ j=1,2,\dots, \ssl$,
$\ell', \ell= 1,2,\dots, \Nop$, be real numbers, and let
$\displaystyle c_i^{[\ell',\ell]}=\sum\limits_{j=1}^{\ssl}a_{ij}^{[\ell',\ell]}$. One
step of a GARK method with an $\Nop$-additive splitting of the
right-hand side of \cref{eq:N-ODE} with $\Nop$ stages takes the form
\begin{align*}
	\mathbf{Y}^{[\ell']}_i & = \mathbf{y}_n +\Delta t
	\sum\limits_{\ell=1}^{\Nop}
	\sum_{j=1}^{\ssl}
	a_{ij}^{[\ell',\ell]}\Fl{\ell}(t_n+c_j^{[\ell',\ell]}\Delta
	t, \mathbf{Y}^{[\ell]}_j),
	\enskip i=1,2,
	\dots,\sslp,\ \ell'=1,2,\dots,\Nop. \\
	\mathbf{y}_{n+1} & = \mathbf{y}_n +\Delta t\sum_{\ell=1}^\Nop \sum_{i=1}^{\ssl} b_i^{[\ell]}\Fl{\ell}(t_n+c_i^{[\ell,\ell]}\Delta t, \mathbf{Y}^{[\ell]}_i).  \\
\end{align*}
The corresponding generalized Butcher tableau is
\begin{equation}\label{FSRKgark_tab}
	\begin{array}{cccc|cccc}
		\cc{1,1}  & \cc{1,2}    &  \cdots  & \cc{1,N}     & 
		\AAA{1,1} & 	\AAA{1,2} & \cdots & 	\AAA{1,N} \\
		\cc{2,1}  & \cc{2,2}    &  \cdots  & \cc{2,N}       & 
		\AAA{2,1} & 	\AAA{2,2} & \cdots & 	\AAA{2,N} \\
		\vdots  & \vdots & \ddots & \vdots    & 
		\vdots  & \vdots & \ddots & \vdots \\
		\cc{\Nop,1}  & \cc{\Nop,2}    &  \cdots  & \cc{\Nop,N}       & 
		\AAA{\Nop,1} & 	\AAA{\Nop,2} & \cdots & 	\AAA{\Nop,N} \\
		\hline 
		&    &    &     & 
		\bb{1}  & \bb{2} & \cdots & \bb{N} \\
	\end{array}
\end{equation}

\begin{remark}
  The diagonal blocks $\AAA{\ell,\ell}$ correspond to the coefficients
  of the Runge--Kutta method for operator $\Fl{\ell}$. The
  off-diagonal entries represent the coupling between operators
  $\Fl{\ell}$ and $\Fl{\ell'}$. As shown in \cite{sandu2015}, GARK
  methods are equivalent to ARK methods; i.e., the Butcher tableau
  \eqref{FSRKgark_tab} is equivalent to a tableau of the form
  \eqref{FSRKeq:arktab}. As shown in \cite{spiteri_wei_FSRK}, an FRSK
  method can be written as an ARK or GARK method. To implement a new
  GARK method in \pythOS, it is sufficient to supply the individual
  blocks $\AAA{\ell,\ell'}$ and $\bb{\ell}$, The blocks
  $\cc{\ell, \ell'}$ are derived from $\AAA{\ell,\ell'}$ such that the
  method is internally consistent, i.e.,
  $\cc{\ell, \ell'} = \AAA{\ell, \ell'}\one^{[s]}$.
\end{remark}

\subsection{Multi-rate methods}
\label{subsec:MR}

Multi-rate methods are specially designed for problems where different
physical components evolve on different time scales, such as problems
seen in numerical weather prediction. The majority of multi-rate
methods \cite{Sandu2019} focus on splitting \cref{eq:N-ODE} into two
operators: a slow ($S$) and a fast ($F$) process, 
\begin{equation} \label{eq:mri2}
	\frac{d\yy}{dt} = \Fl{S}(t, \yy) + \Fl{F}(t, \yy). 
\end{equation}
The slow process $\Fl{S}$ is solved with a macro time step $\Dt$, and
the fast process $\Fl{F}$ is solved with a smaller micro time step
$\displaystyle \frac{\Dt}{M}$, where $M$ is the number of micro-steps
per macro-step. As shown in \cite{gunther_multirate_2016}, if
Runge--Kutta methods are used as the integration methods for the slow
and fast processes, the corresponding multi-rate method can be written
as a GARK method. We refer to this method as a multi-rate GARK (MrGARK)
method. Let $\{ \AAA{S,S}, \bb{S,S}, \cc{S,S} \}$ and
$\{ \AAA{F,F}, \bb{F,F}, \cc{F,F} \}$ be the coefficients of the
$s^{[S]}$- and $s^{[F]}$-stage Runge--Kutta methods corresponding to
the slow and fast processes, respectively, and let
$\{ \AAA{F,S,\faststepind}, \cc{F,S,\faststepind}\}$ and $\{ \AAA{S,F,\faststepind}, \cc{S,F,\faststepind} \}$,
for $\faststepind =1,2,\dots, M$, be the coupling coefficients between the
fast-slow and slow-fast processes. An MrGARK method
applied to \cref{eq:mri2} that advances the solution from $t_n$ to
$t_{n+1} = t_n+ \Dt$ takes the form

\begin{align*}
		\YY^{[S]}_i & = \yy_n +\Dt \sum\limits_{j=1}^{s^{[S]}}a_{i,j}^{[S,S]}\Fl{S}\left(t_n+c_j^{[S,S]}\Dt, \YY_j^{[S]} \right) \\
		& \quad + \frac{\Dt}{M}\sum\limits_{ \faststepind =1}^{M}\sum\limits_{j=1}^{s^{[F]}}a_{i,j}^{[S,F, \faststepind]}\Fl{F}\left(t_n+\left(c_j^{[S,F, \faststepind]}+ \faststepind\right)\frac{\Dt}{M}, \YY_j^{[F, \faststepind]}\right),
		\enskip i=1,2,
		\dots,s^{[S]},
\end{align*}
where
\begin{align*}
		\YY^{[F, \faststepind]}_i & = 
		\yy_n + \Dt \sum\limits_{l=1}^{ \faststepind-1} \sum\limits_{j=1}^{s^{[F]}} b^{[F,F]}_j \Fl{F}\left(t_n+\left(c_i^{[F,F]}+ \faststepind \right)\frac{\Dt}{M}, \YY_i^{[F, \faststepind]}\right)\\
		& \quad + \Dt \sum\limits_{j=1}^{s^{[S]}}a_{i,j}^{[F,S, \faststepind]}\Fl{S}\left(t_n+c_j^{[F,S, \faststepind]}\Dt,\YY_j^{[S]}\right) \\
		& \quad + \frac{\Dt}{M}
		\sum\limits_{j=1}^{s^{[F]}}a_{i,j}^{[F,F]}\Fl{F}\left(t_n+\left(c_j^{[F,F]}+ \faststepind\right)\frac{\Dt}{M}, \YY_j^{[F, \faststepind]}\right), \enskip i=1,2,\dots,s^{[F]},  \faststepind=1,2,\dots,M,
\end{align*}

\begin{align*}
	\yy_{n+1} & = \yy_n + \frac{\Dt}{M} \sum\limits_{ \faststepind = 1}^M \sum\limits_{i=1}^{s^{[F]}} b_i^{[F,F]} \Fl{F}\left(t_n+\left(c_i^{[F,F]}+ \faststepind\right)\frac{\Dt}{M}, \YY_i^{[F, \faststepind]}\right) \\ 
	& \quad +  \Dt\sum\limits_{i=1}^{s^{[S]}}b_i^{[S,S]}\Fl{S}\left(t_n+c_i^{[S,S]}\Dt, \YY_i^{[S]}\right). 
\end{align*}

The corresponding generalized Butcher tableau is

\begingroup
\renewcommand\arraystretch{1.75}
\begin{equation}\label{FSRKmr_tab}
  \begin{array}{cccc|c}
		\frac{1}{M}\AAA{F,F} & 	0 & \cdots & 	0 & \AAA{F,S,1}\\
		\frac{1}{M}\one\bb{F,F}\,^T & \frac{1}{M}\AAA{F,F} & \cdots & 0 & \AAA{F,S,2} \\
		\vdots  & \vdots & \ddots & \vdots & \vdots\\
		\frac{1}{M}\one\bb{F,F}\,^T &\frac{1}{M}\one\bb{F,F}\,^T & \cdots & 	\frac{1}{M}\AAA{F,F} & \AAA{F,S,M} \\
    \hline
    \frac{1}{M}\AAA{S,F,1} & \frac{1}{M}\AAA{S,F,2} & \cdots & \frac{1}{M}\AAA{S,F,M} & \AAA{S,S} \\
    \hline
    \frac{1}{M}\bb{F,F}\,^T & \frac{1}{M}\bb{F,F}\,^T & \cdots & \frac{1}{M}\bb{F,F}\,^T & \bb{S,S}\,^T
  \end{array}
\end{equation}
\endgroup

If the fast process is solved exactly (in principle ``with an
infinitesimally small step size''), the resulting method is called a
multi-rate infinitesimal method.  In \pythOS, we have
implemented MRI-GARK methods \cite{Sandu2019} to solve
\cref{eq:mri2}. An MRI-GARK method is based on an $\mrislowstage$-stage
``slow'' Runge--Kutta method with non-decreasing abscissae
$c^{S}_1 \leq c^{S}_2 \leq \cdots \leq c^{S}_{\mrislowstage}
\leq 1$ and lower-triangular coupling matrices $\GGamma{k}$. An
MRI-GARK method applied to \cref{eq:mri2} advances the solution from
$t_n$ to $t_{n+1} = t_n+ \Dt$ as shown in \cref{alg:mri2}.

\begin{algorithm}[!hbtp]
  \caption{\label{alg:mri2} Algorithm for taking a single step for an 
    MRI-GARK method.}
	
	\begin{algorithmic}[1]

		\State $\YY_1^{[S]} = \yy_n$
		
		\For {$i=2$ to $\mrislowstage$}
			\State $v(0) := \YY_{i-1}^{[S]}$  and $T_{i-1}:= t_n+c_{i-1}^{[S]}\Dtn$,
			\State $\Delta c_{i}^{[S]} = c_i^{[S]} - c_{i-1}^{[S]}$
			\State Solve  $v'(\theta) = \Delta c_i^{[S]} \Fl{F}(T_{i-1} + \Delta c_i^{[S]}\theta, v(\theta)) + g(\theta)$ for $\theta \in [0,\Dtn]$,
			\State where $g(\theta) = \sum\limits_{j=1}^i \gamma_{i,j}(\frac{\theta}{\Dt}) \Fl{S}(t_n+c_j^{[S]} \Dtn, \YY_j^{[S]})$,
			\State $\YY_i^{[S]} := v(\Dtn)$,
		\EndFor

		\State \textbf{Return} $\yy_{n+1}:= \YY_{\mrislowstage}^{[S]}$
	\end{algorithmic}
\end{algorithm}
Here,
$\gamma_{i,j}(\tau) := \sum\limits_{k\geq 0} \gamma_{i,j}^{\{k\}}
\tau^k$, where $\gamma_{i,j}^{\{k\}}$ are entries of
$\GGamma{k}$. The matrices $\GGamma{k}$ are derived from order
conditions for MRI-GARK methods and determine the coupling between
the slow and fast components.

Some methods \cite{Chinomona2021} split the slow process further into
stiff and non-stiff operators. The stiff operator is usually solved
implicitly ($I$) and the non-stiff operator is solved explicitly ($E$);
i.e.,
\begin{equation}\label{eq:mri3}
	\frac{d\yy}{dt} = \Fl{I}(t, \yy) + \Fl{E}(t, \yy) + \Fl{F}(t, \yy).
\end{equation}
To solve \cref{eq:mri3}, the IMEX-MRI-GARK methods proposed in
\cite{Chinomona2021} are implemented in \pythOS. The IMEX-MRI-GARK
methods are an extension of the MRI-GARK methods in \cite{Sandu2019}. In
addition to an $\mrislowstage$-stage "slow" Runge--Kutta method with
non-decreasing abscissae
$c_1^{[S]} \leq c_2^{[S]} \leq \cdots \leq c_{\mrislowstage}^{[S]}
\leq 1$, the IMEX-MRI-GARK methods have two sets of coupling matrices
$\GGamma{k}$ and $\OOmega{k}$. Here, the $\GGamma{k}$ are lower
triangular, and the $\OOmega{k}$ are strictly lower triangular. Let
$\gamma_{i,j}(\tau) := \sum\limits_{k\geq 0} \gamma_{i,j}^{\{k\}}
\tau^k$ and
$\omega_{i,j}(\tau) := \sum\limits_{k\geq 0} \omega_{i,j}^{\{k\}}
\tau^k$, where $\gamma_{i,j}^{\{k\}}$ and $\omega_{i,j}^{\{k\}}$ are the
entries of $\GGamma{k}$ and $\OOmega{k}$, respectively. An
IMEX-MRI-GARK method applied to \cref{eq:mri3} advances the solution
from $t_n$ to $t_{n+1} = t_n+ \Dtn$ as shown in \cref{alg:mri3}.

\begin{algorithm}[!hbtp]
	\caption{\label{alg:mri3} Algorithm for taking a single
		step for an IMEX-MRI-GARK method.}
	
	\begin{algorithmic}[1]

		\State $\YY_1^{[S]} = \yy_n$
		
		\For {$i=2$ to $\mrislowstage$}
		\State $v(0) := \YY_{i-1}^{[S]}$  and $T_{i-1}:= t_n+c_{i-1}^{[S]}\Dtn$,
		\State $\Delta c_{i}^{[S]} = c_i^{[S]} - c_{i-1}^{[S]}$
		\State Solve  $v'(\theta) = \Delta c_i^{[S]} \Fl{F}(T_{i-1} + \Delta c_i^{[S]}\theta, v(\theta)) + g(\theta)$ for $\theta \in [0,\Dtn]$,
		\State where $g(\theta) = \sum\limits_{j=1}^i \gamma_{i,j}(\frac{\theta}{\Dtn}) \Fl{I}(t_n+c_j^{[S]} \Dtn, \YY_j^{[S]}) + \sum\limits_{j=1}^{i-1} \omega_{i,j}(\frac{\theta}{\Dtn}) \Fl{E}(t_n+c_j^{[S]} \Dtn, \YY_j^{[S]})$,
		\State $\YY_i^{[S]} := v(\Dtn)$,
		\EndFor
		
		\State \textbf{Return} $\yy_{n+1}:= \YY_{\mrislowstage}^{[S]}$
	\end{algorithmic}
\end{algorithm}
    
\section{Structure of \pythOS}
\label{sec:structure}
The \pythOS\ library contains five separate solvers that involve
splitting.  Each solver may be used by importing the corresponding
main file: {\tt fractional\_step} for fractional-step methods, {\tt
  additive\_rk} for additive Runge--Kutta methods, {\tt gark\_methods}
for generalized additive Runge--Kutta methods, {\tt multirate}
for multi-rate methods and {\tt
  multirate\_infinitesimal} for multi-rate infinitesimal methods.  The
libraray is available on GitHub at \url{https://github.com/uofs-simlab/pythOS}.

%

Examples of the basic form of the call to the solvers are provided in
each of the following sections, and a complete description of the
optional arguments is provided in the documentation accompanying the
software. Unless otherwise noted, in all subsections of this chapter,
the ODEs are solved with initial condition {\tt y\_0} over the time
interval [{\tt t0}, {\tt tf}] with step size {\tt dt}.
  
\subsection{Fractional-step methods}

The fractional-step solver can be imported using

{\tt import fractional\_step as fs},\\
which provides the solver function {\tt fractional\_step}.  This
function takes as inputs a list of operators, an initial condition,
initial and final times, a time step, a string (or list of lists) for
the splitting method, and a dictionary defining the sub-integrators.
It also optionally takes inputs to save intermediate results, control
the behaviour of the sub-integrators, perform dynamic linearization,
or control splitting error.

For example, using functions {\tt f1} and {\tt f2} for the operators,
the call to the fractional-step solver using Godunov splitting,
forward Euler ({\tt FE}) to solve the first operator, and backward
Euler ({\tt BE}) to solve the second operator is as follows:

\begin{lstlisting}
result = fs.fractional_step([f1, f2], dt, y_0, t0, tf, 
"Godunov", methods={(1,): "FE", (2,): "BE"})
\end{lstlisting}

The \pythOS\ library provides options to use Runge--Kutta, exponential
propagation iterative (EPI) \cite{Tokman2006}, adaptive, or analytic
sub-integrators.  Each choice is specified by changing the
corresponding entry in the {\tt methods} dictionary.  There are
Runge--Kutta methods included in \pythOS, and there is the option to
use methods from \Irksome\ for problems using finite element
discretizations. For EPI methods, the underlying solver may be
\texttt{kiops} \cite{Gaudreault2018}
and can be specified on a per operator basis.  For an adaptive
sub-integration, the SUNDIALS solvers and the solvers from {\tt
  scipy.integrate.solve\_ivp} are built-in \pythOS\ options, as well
as a collection of embedded Runge--Kutta methods.  The choice of
adaptive solver can be specified on a per operator basis as well.
Alternatively, we can specify embedded Runge--Kutta methods in {\tt
  ivp\_methods} to solve IVPs with complex-valued time steps. It is
worth noting, however, that {\tt scipy.integrate.solve\_ivp} does not
support complex-valued time for IVPs.

A fractional-step method is defined by a list of lists that defines
the size of the steps.  Entry $i$ of list $j$ defines the fractional
step for step $j$ of operator $i$.  Before solving, the software
combines this information with the method specification and the
operator specification to create a list of non-trivial steps to take.
A new method may be specified by providing an appropriate list of
lists.

\subsection{Additive Runge--Kutta methods}

The additive Runge--Kutta solver in \pythOS\ can be imported using

{\tt import additive\_rk as ark}\\
which provides the solver function {\tt ark\_solve}.  This function
takes as inputs a list of operators, initial condition, initial and
final times, a time step, and a list of Butcher tableaux that define
the method.  It also takes optional arguments that can be used to
control the non-linear solver used for implicit
stages, set tolerances for adaptive step methods, save intermediate
output, specify boundary conditions if using finite element
capabilities from \Firedrake, and perform dynamic linearization if
desired.

For example, using functions {\tt f1} and {\tt f2} for the operators,
and {\tt tableau1} and {\tt tableau2} as the Butcher tableaux
corresponding to each operator, the call to the additive Runge--Kutta
solver is as follows:
\begin{lstlisting}
result = ark.ark_solve([f1, f2], dt, y_0, t0, tf, 
[tableau1, tableau2])
\end{lstlisting}


The additive Runge--Kutta method is defined by a list of Butcher
tableaux, each of type {\tt Tableau}.  These may easily be created by
specifiying the arrays $a$, $b$, and $c$ that define the tableau.  An
embedded additive Runge--Kutta method is defined by all the Butcher
tableaux being of type {\tt EmbeddedTableau}, which are created by
specifying the arrays $c$, $a$, $b$, $\tilde{b}$, and the order value
$p$.  The solver requires the size of each tableau to be the same for
each operator.  The solver takes the structure of the underlying
tableaux into account in order to apply an optimized implementation
for diagonally implicit tableaux.

\subsection{Generalized additive Runge--Kutta methods}

The generalized additive Runge--Kutta method solver in \pythOS\ can be
imported as

{\tt import gark\_methods as gark}\\
which provides the
solver function {\tt gark\_solve}.  This function takes as inputs a
list of operators, initial condition, initial and final times, a time
step, a list of lists of \numpy\ arrays that define each section
$\AAA{i,j}$ of the GARK tableau, and a list of \numpy\ arrays that
defines the arrays $\bb{i}$ for the method.  It also takes optional
arguments that can be used to control the solver used for implicit
stages, save intermediate output, specify boundary conditions if using
finite element capabilities from \Firedrake, and perform dynamic
linearization, if desired.

For example, using functions {\tt f1} and {\tt f2} for the operators,
and {\tt A11}, {\tt A12}, {\tt A21}, {\tt A22} as the arrays
$\AAA{i,j}$, and {\tt b1} and {\tt b2} as the arrays $\bb{i}$, the
call to the generalized additive Runge--Kutta solver is as follows:
\begin{lstlisting}
result = gark.gark_solve([f1, f2], dt, y_0, t0, tf,
[[A11, A12], [A21, A22]], [b1, b2])
\end{lstlisting}

To solve this type of problem, the input arrays are converted into the
corresponding additive Runge--Kutta tableaux and solved using the
additive Runge--Kutta solver.

\subsection{Multi-rate methods}

\subsubsection{Multi-rate GARK methods}{\ }

The multi-rate GARK solver can be imported using 

{\tt import multirate as mr}, \\ 
which provides the solver function {\tt multirate\_solve}.
This function takes as inputs a slow operator, a fast operator, a time step,
an initial condition, initial and final times, and method information.  It also
takes optional arguments that can be used to control the solver used for the
implicit stages, the solver used for the fast stages, save
intermediate output, and specify boundary conditions if using finite
element capabilities from \Firedrake.

For example, using functions {\tt fS} and {\tt fF} for the slow and
fast operators, respectively, and the MrGARK-EX2-IM2 method
\cite{sarshar2019} with $M$ fast steps per slow time step {\tt dt},
the call to the multi-rate solver is as follows:
\begin{lstlisting}
result = mr.multirate_solve(y_0, t0, dt, tf, 
mr.mrgark_ex2_im2, M, fS, fF)
\end{lstlisting}

Other MrGARK methods may be used by creating a {\tt Multirate} object.  The
object contains the \numpy\
arrays {\tt  A\_FF}, {\tt A\_SS}, {\tt b\_F}, and {\tt b\_S} as matrices
$\AAA{F,F}$, $\AAA{S,S}$, $\bb{F,F}$, and $\bb{S,S}$ and functions {\tt
  A\_FS} and {\tt A\_SF} in terms of $\faststepind$ and $M$ to compute matrices
$\AAA{F,S,\faststepind}$ and $\AAA{S,F,\faststepind}$.

To solve this type of problem, the input arrays are converted into the
corresponding generalized additive Runge--Kutta tableaux, reordered if
the method is decoupled, and solved using the
generalized additive Runge--Kutta solver.

\subsubsection{Multi-rate infinitesimal methods}{\ }

The multi-rate infinitesimal solver can be imported using

{\tt import
  multirate\_infinitesimal as mri}\\
which provides the solver function {\tt
  multirate\_infinitesimal\_solve}.  This function takes as inputs an initial
condition, initial and final times, a time step, the method (of type\\
{\tt Multirate\_Infinitesimal}), an implicit slow operator, a fast operator, and
optionally an explicit slow operator.  It also takes optional
arguments that can be used to control the solver used for the
implicit stages, control the solver used for the fast stages, save
intermediate output, and specify boundary conditions if using finite
element capabilities from \Firedrake.

For example, using functions {\tt fS} and {\tt fF} for the slow and fast
operator, and
{\tt mri\_sdirk3} as the method (which is provided in the
{\tt multirate\_infinitesimal}
file), the call to the multi-rate solver is as follows:
\begin{lstlisting}
result = mri.multirate_infinitesimal_solve(y_0, t0, dt, tf, 
mri.mri_sdirk3, fS, fF)
\end{lstlisting}
  
The solvers from {\tt scipy.integrate.solve\_ivp}, the non-split SUNDIALS
solvers, and the embedded Runge--Kutta methods built into \pythOS\ may all
be used for solving the fast scale.
New methods may be used by creating a {\tt Multirate\_Infinitesimal} object.
If the new method is an MRI-IMEX method, a new
{\tt Multirate\_Infinitesimal} object is created
by providing the arrays $c_i^{[S]}$, $\gamma_{ij}^{\{k\}}$, and
$\omega_{ij}^{\{k\}}$.  If the new method is not split at the slow scale,
the {\tt Multirate\_Infinitesimal} object is created by providing the arrays $c_i^{[S]}$ and
$\gamma_{ij}^{\{k\}}$.  We note that the solver requires that there are
no implicit solves on stages with adaptive integrations, so
$\Delta c_i^{[S]} =0$ if $\omega_{ij}^{\{k\}} \neq 0$, $j \geq i$ or
$\gamma_{ij}^{\{k\}} \neq 0$, $j \geq i$.

\section{Examples}
\label{sec:examples}

In this section, we demonstrate the functionality of the main \pythOS\
solvers while illustrating some perhaps under-appreciated aspects of
operator-splitting methods, such has high convergence order (despite
backward steps) and the use of methods with complex coefficients,
especially for the direct solution of complex ODEs.


\subsection{2D advection-diffusion-reaction problem}

We first consider a two-dimensional advection-diffusion-reaction
problem,

\begin{equation} \label{eq:ADR} u_t = - \alpha \left(\nabla \cdot u\right) + \epsilon
  \left(\nabla^2 u\right) + \gamma u \left(u-1/2\right)\left(1-u\right),
\end{equation}
with homogeneous Neumann boundary conditions and initial conditions,
\begin{equation*}
  u\left(x, y, 0\right) = 256\left(xy\left(1-x\right)\left(1-y\right)\right)^2 + 0.3,
\end{equation*}
and solved on $t \in [0, 0.1]$, $(x, y) \in [0, 1]^2$ with parameters
$\alpha = -10$, $\epsilon=1/100$, $\gamma = 100$.  For the spatial
discretization, we use a uniform grid with grid size
$\displaystyle \Dx = \Dy = 1/40$. Spatial derivatives are
discretized using the continuous Galerkin finite element method with
quadratic elements as implemented in \Firedrake.

\begin{figure}[htbp]
  \centering \includegraphics[width =
  \textwidth]{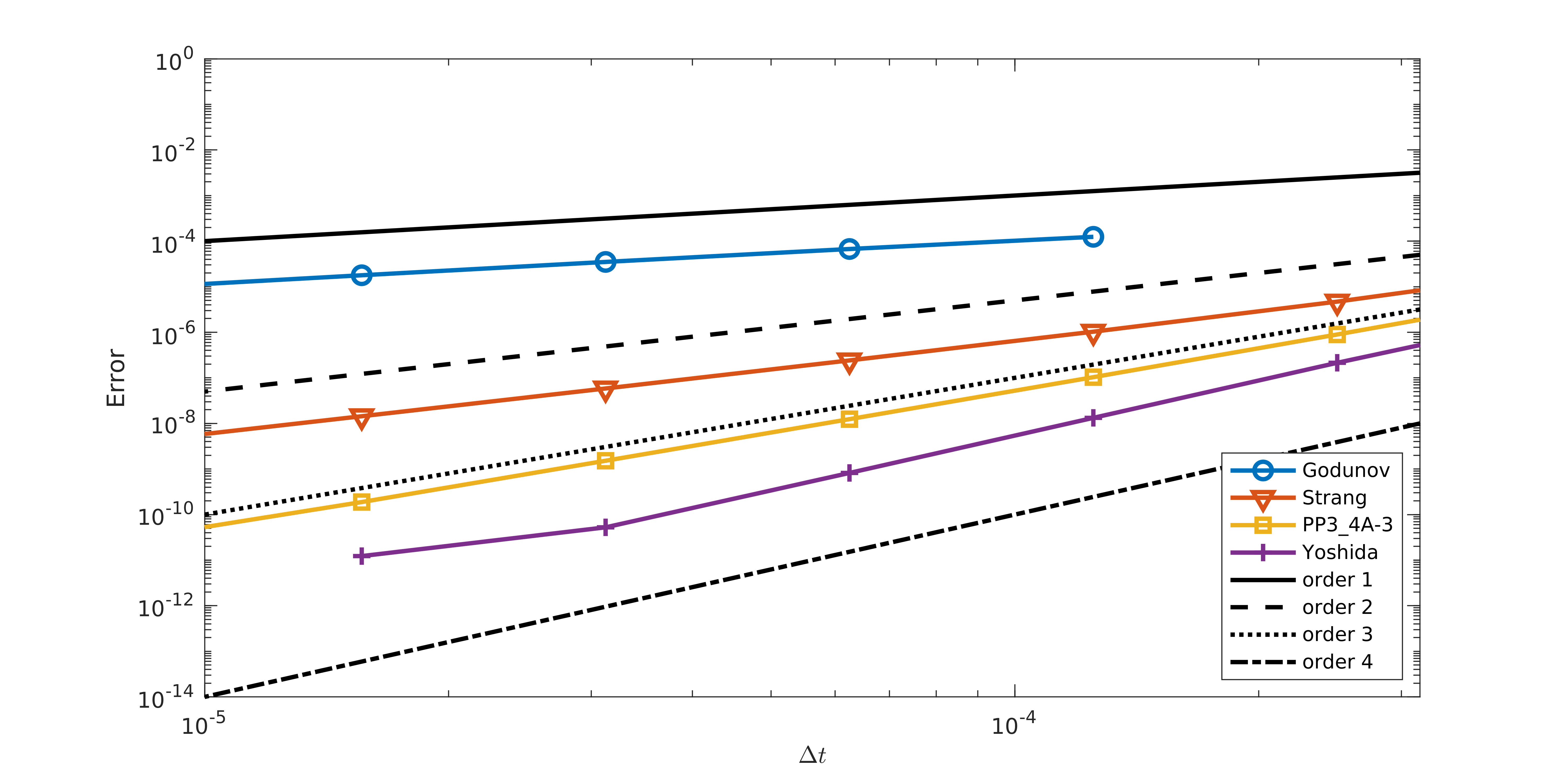}
  \caption{Convergence of the Godunov, Strang, PP3\_4A-3, and Yoshida
    splitting methods applied to the 2D advection-diffusion-reaction
    problem \cref{eq:ADR}}
	\label{fig:ADR_2D_OS_RK4_convergence}
\end{figure}

We solve this problem by first splitting \cref{eq:ADR} into three
operators: advection, diffusion, and reaction. In particular,
\begin{align*} 
&&\Fl{1}(u) = -\alpha\left(\nabla \cdot u\right), && \Fl{2}(u) =  \epsilon
\left(\nabla^2 u\right), && \Fl{3}(u) = \gamma u \left(u-1/2\right)(1-u). 
\end{align*}
The following $3$-split operator-splitting methods were tested: the
first-order Godunov method, the second-order Strang method, a
third-order palindromic method PP3\_4A-3 \cite{auzinger2017}, and the
fourth-order Yoshida method~\cite{Yoshida1990}.  The coefficients for
these methods are given in \cref{sec:os_methods_coeff}.  The
sub-systems are solved using sub-integration methods of the same order
of accuracy as the splitting method, namely, forward Euler for the
Godunov splitting, the two-stage, second-order explicit Heun method
for the Strang splitting, the explicit third-order Kutta method
\cite{kutta1901} for the PP3\_4A-3 splitting, and the classical
four-stage, fourth-order explicit Runge--Kutta method for the Yoshida
splitting.  The reference solution is obtained using the
Dormand--Prince~5(4) adaptive solver from \pythOS\ with \rtol = \atol
= \num{1e-14}. The error is measured using the $\ell_2$-norm of the
difference between the numerical solution and the reference solution
at $\tf = 0.1$. \Cref{fig:ADR_2D_OS_RK4_convergence} shows that all
four methods exhibit the expected order of convergence.

\begin{figure}[htbp]
  \centering \includegraphics[width =
  \textwidth]{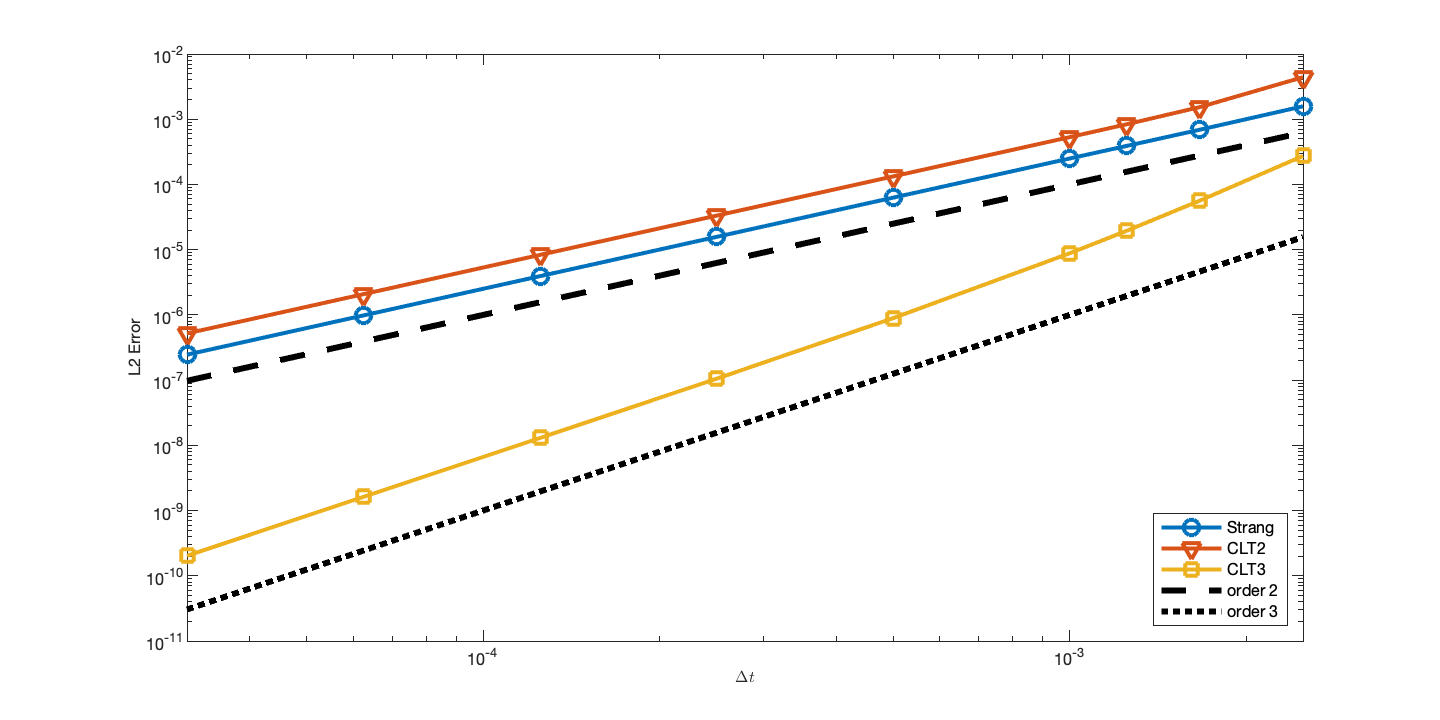}
  \caption{Convergence of the Strang, CLT2, and CLT3 splitting methods
    applied to the 4-split 2D advection-diffusion-reaction problem
    \cref{eq:ADR}}
	\label{fig:adr_2d_4split_convergence}
\end{figure}

To demonstrate the ability to solve an ODE using an arbitrary number
of operators, we also present the result of solving \cref{eq:ADR} as
$4$-split problem with the following four
operators:
\begin{align*}
	 \Gl{1}(u) = - \alpha\, (\nabla \cdot u),  
		&& \Gl{2}(u)  = \epsilon u_{xx},  
		&& \Gl{3}(u)  = \epsilon u_{yy},  
		&& \Gl{4}(u)  = \gamma u \left(u-1/2\right)\left(1-u\right).
\end{align*}
The $4$-split problem is solved using the second-order Strang method,
the second-order complex-valued \CLTTwo\ (CLT2) method, and the
third-order complex-valued \CLTThree\ (CLT3) method
\cite{spiteri_wei_clt}. The coefficients of these methods can be found
in \cref{sec:os_methods_coeff}. In \pythOS, these methods can be
called without specifying the number of operators $\Nop$. The code
identifies the number of operators and executes the corresponding
method. \Cref{fig:adr_2d_4split_convergence} shows that all three
methods exhibit the expected order of convergence.

\subsection{Complex-valued ODE}

We consider a complex-valued ODE,

\begin{equation}
  \label{complex_ode}
  \dv{u}{t} = iu + 0.1u - 0.1u^3,
\end{equation}
with initial condition $u(0) = 0.1$ and solved on $t\in [0, 100]$.
\Cref{complex_ode} can be solved as a complex-valued differential
equation or as a system of real-valued differential equations.

To solve \cref{complex_ode} as a complex-valued ODE using $3$-split
operator-splitting methods, we split \cref{complex_ode} as
\begin{align*} \label{complex_ode_split}
	&& \complexde^{[1]} (u) = iu, && \complexde^{[2]}(u) = 0.1 u, && \complexde^{[3]}(u) = -0.1u^3.
\end{align*}

To solve \cref{complex_ode} as a system of real-valued differential
equations, let $u = x+iy$. Then \cref{complex_ode} can be rewritten as
a system of real-valued differential equations in $x$ and $y$:
\begin{equation} \label{complex_ode_real}
	\left\{ 
	\begin{aligned}
		\dv{x}{t} & = - y+  0.1x  +0.3xy^2 -0.1x^3,  \\
		\dv{y}{t} & = \hspace*{1.5ex}x+0.1y -0.3x^2y + 0.1y^3.
	\end{aligned}  \right.  
\end{equation}
To obtain a reference solution to \cref{complex_ode}, the system is
converted to \cref{complex_ode_real} and solved using MATLAB's {\tt
  ode45} with {\rtol} = \num{1e-16} and {\atol} = \num{1e-15}. The
error is calculated using the mixed root-mean-squre (MRMS)
error~\cite{Marsh2012} defined by
$$ \text{MRMS} = \sqrt{\frac{1}{n_t}\sum_{i=1}^{n_t}\left(\frac{u_i-u^{\text{ref}}_i}{1+\vert
      u^{\text{ref}}_i\vert}\right)^2},$$ where $u_i$,
$u^{\text{ref}}_i$ denote the numerical and reference solution,
respectively, at time $t_i$ and $n_t$ denotes the total number of time
points at which the error is calculated. For this example, $t_i =
i$, $i=1,2,\dots,100$.

\Cref{complex_ode_real} can be split into three operators, 
\begin{align*}
&&	\realde^{[1]}(x,y) = \begin{bmatrix}
		-y \\
		\hspace*{1ex}x \\
	\end{bmatrix},  
&& \realde^{[2]}(x,y) = \begin{bmatrix}
	0.1x \\ 
	0.1y \\
\end{bmatrix},
&& 
	\realde^{[3]}(x,y) = \begin{bmatrix}
	\hspace*{1.75ex}0.3xy^2-0.1x^3 \\
	-0.3x^2y + 0.1y^3
\end{bmatrix}.
\end{align*}
We apply the second-order Strang method, second-order complex-valued
CLT2 method, and third-order complex-valued CLT3 method using
$\complexde^{[1]}$, $\complexde^{[2]}$, $\complexde^{[3]}$, and
$\realde^{[1]}$, $\realde^{[2]}$, $\realde^{[3]}$. The coefficients of
these splitting methods can be found in \cref{sec:os_methods_coeff}.

\begin{figure}[htbp]
	\centering
	\includegraphics[width = \textwidth]{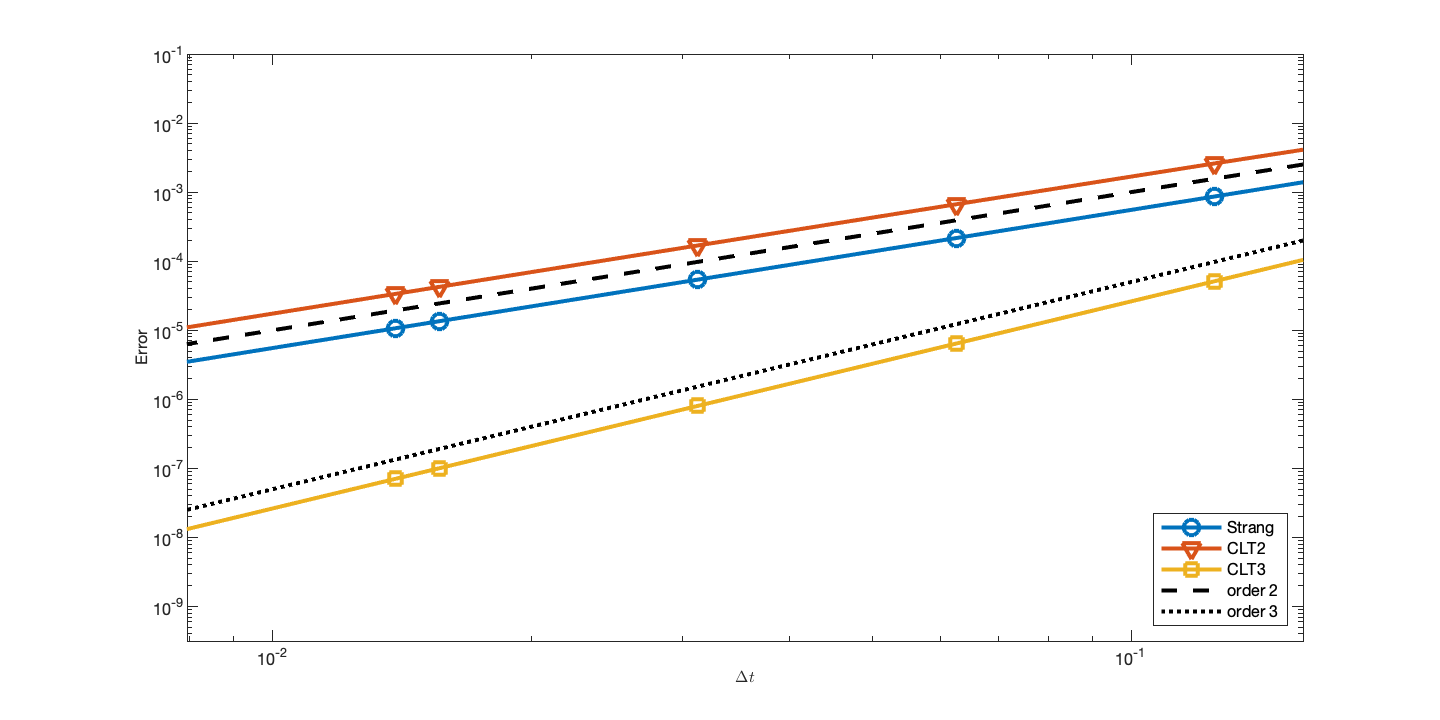}
	\caption{Convergence of the 3-split versions of Strang,
      CLT2, and CLT3 splitting applied to the complex ODE
      \cref{complex_ode}.}
	\label{fig:complex_ode_convergence}
\end{figure}

\begin{figure}[htbp]
	\centering
	\includegraphics[width = \textwidth]{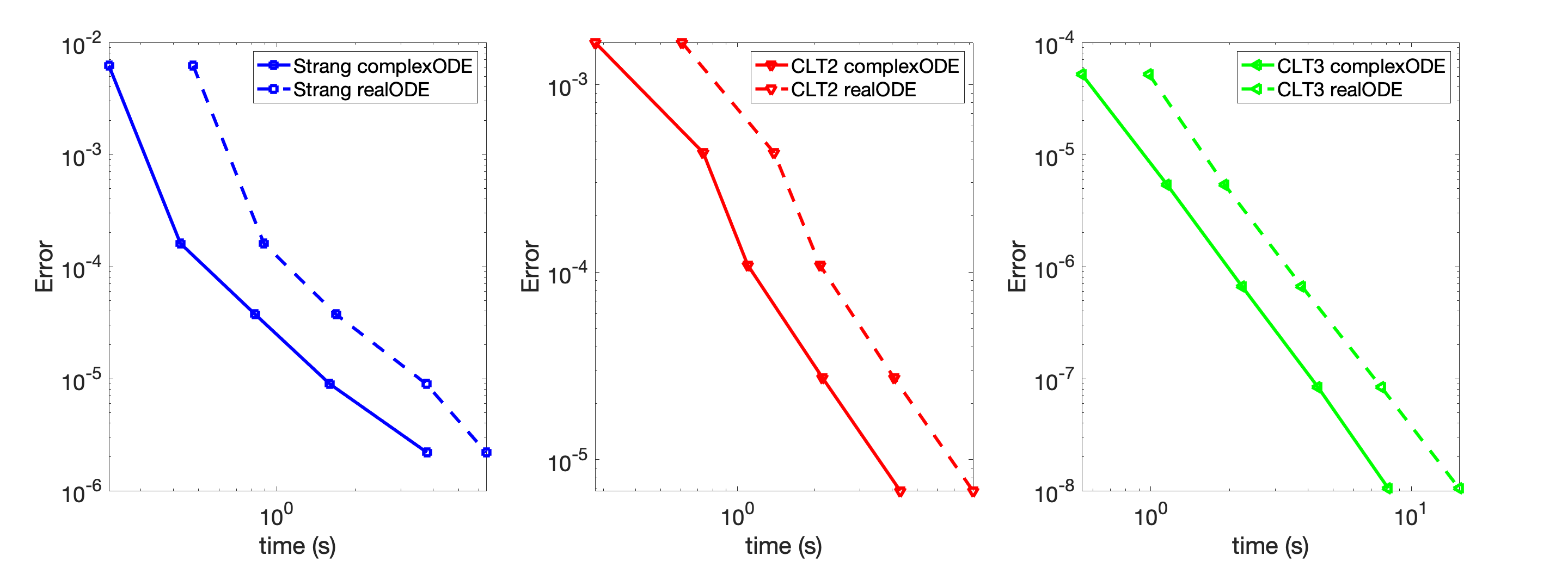}
	\caption{Work-precision diagrams of Strang, CLT2, and CLT3
      with RK3 sub-integrators applied to the complex ODE
      \cref{complex_ode} and the system of real ODEs
      \cref{complex_ode_real}.}
	\label{fig:complex_ode_work_precision_rk3}
\end{figure}

In order to investigate the effect of the choice of sub-integration
method on performance, two experiments are conducted: the sub-systems
$\displaystyle \dv{u^{[i]}}{t} =
\complexde^{[i]}\left(u\right)$ and $\begin{bmatrix} \displaystyle \dv{x^{[i]}}{t} \\ \displaystyle \dv{y^{[i]}}{t} \\
\end{bmatrix} = \realde^{[i]}\left(x,y\right)$ for $i=1,2,3$ are
solved using the Kutta's third-order method (RK3). 
\Cref{fig:complex_ode_convergence} shows that all
three methods exhibit the expected rates of convergence.
Moreover, as shown in the work-precision diagram
\cref{fig:complex_ode_work_precision_rk3}, it is noticeably more
efficient for this problem to directly solve the complex-valued ODE
rather than the real-valued system of ODEs (even with a real
operator-splitting method).

\subsection{Stiff Brusselator problem}

We consider the stiff variant of the one-dimensional Brusselator
problem from~\cite{Chinomona2021}:

\begin{equation}\label{eq:stiff_brusselator}
\begin{aligned}
  & u_t = \alpha_u u_{xx} + \rho_u u_x + a - (w + 1)u + u^2v,\\
  & v_t = \alpha_v v_{xx} + \rho_v v_x + wu - u^2v,\\
  & w_t = \alpha_w w_{xx} + \rho_w w_x + (b-w)/\epsilon - wu,
\end{aligned}
\end{equation}
with stationary boundary conditions
\begin{align*}
  u_t(t, 0) = u_t(t, 1) = v_t(t, 0) = v_t(t, 1) = w_t(t, 0) = w_t(t,1) = 0
\end{align*}
and initial values
\begin{align*}
  u(0,x) &= a + 0.1 \sin(\pi x),\\
  v(0,x) &= b/a + 0.1 \sin(\pi x),\\
  w(0,x) &= b + 0.1 \sin(\pi x),
\end{align*}
solved on $t \in [0, 3]$ and $x \in [0, 1]$ with parameters
$\alpha_j = 10^{-2}$, $\rho_j = 10^{-3}$, $a = 0.6$, $b = 2$, and
$\epsilon = 10^{-3}$.  We discretize using a uniformly distributed spatial
grid of $n_x=201$ grid points, and use central differences to approximate the
spatial derivatives.  A reference solution was generated using {\tt RK45}
from \scipy's {\tt solve\_ivp} with {\rtol} = \num{1e-16} and
{\atol} = \num{1e-16}. 

\begin{figure}[htbp]
  \centering
  \includegraphics[width=5in]{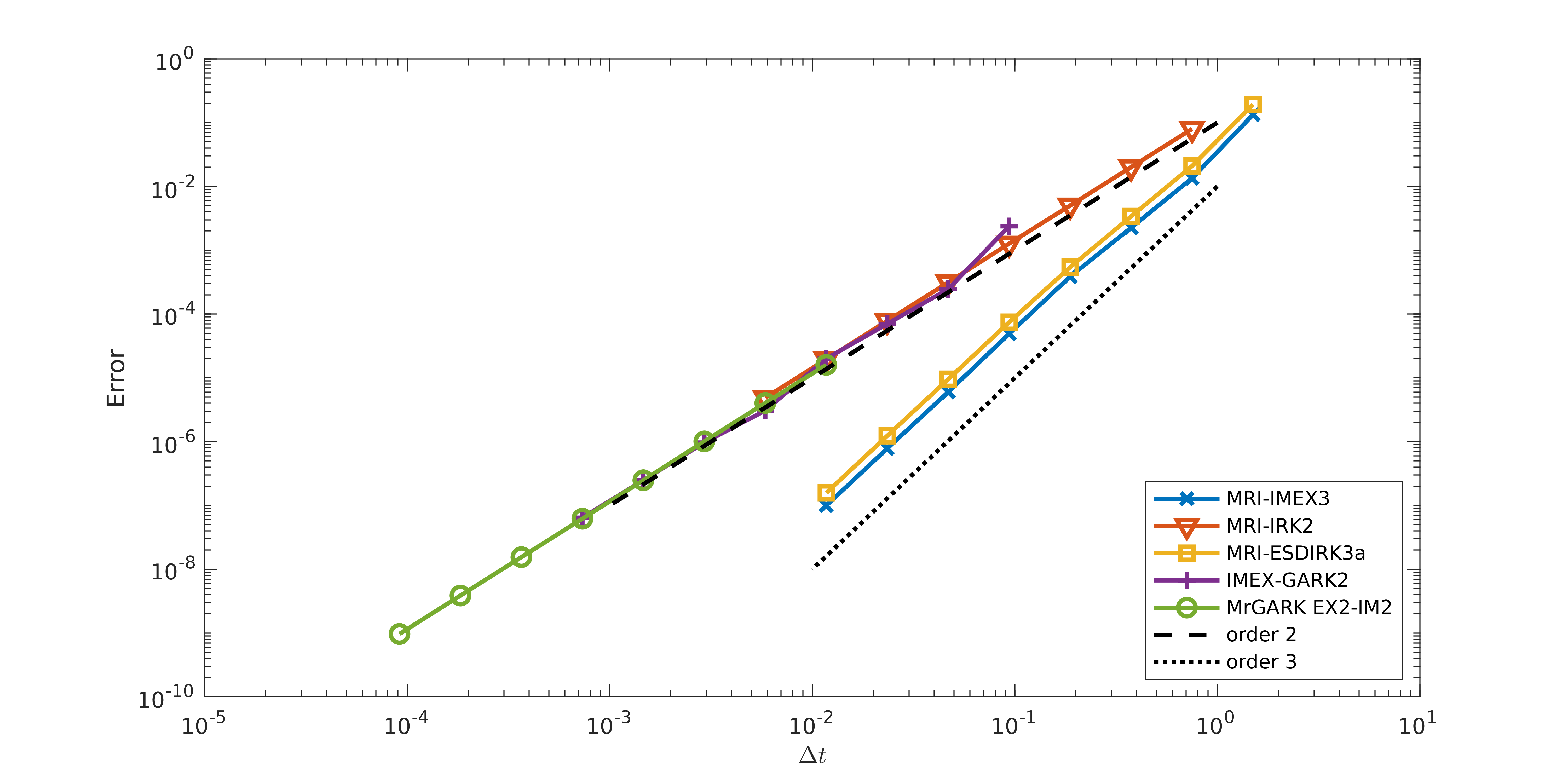}
  \caption{Convergence of multi-rate infinitesimal and generalized additive
    Runge--Kutta methods applied to the stiff Brusselator problem
    \cref{eq:stiff_brusselator}.}
  \label{fig:stiff_brusselator_convergence}
\end{figure}

For multi-rate methods, the right-hand side of
~\cref{eq:stiff_brusselator} is split into the fast and slow
operators:
\begin{equation*} \label{stiff_brus_fs}
\begin{aligned}
	\Fl{S}(u,v,w) = \begin{bmatrix} \rho_u u_x +\alpha_u u_{xx} \\ \rho_v v_x + \alpha_v v_{xx}  \\ \rho_w w_x + \alpha_w w_{xx} \end{bmatrix},
	\Fl{F}(u,v,w) = \begin{bmatrix} a - (w + 1)u + u^2v\\ wu - u^2v\\ (b-w)/\epsilon - wu \end{bmatrix}.
\end{aligned}
\end{equation*}
The slow operator can be further split into the stiff operator $\Fl{I}$ and non-stiff operator $\Fl{E}$: 
\begin{equation*} \label{stiff_brus_fie}
\begin{aligned}
  \Fl{E}(u,v,w) = \begin{bmatrix} \rho_u u_x \\ \rho_v v_x \\ \rho_w w_x \end{bmatrix},
  \Fl{I}(u,v,w) = \begin{bmatrix} \alpha_u u_{xx} \\ \alpha_v v_{xx} \\ \alpha_w w_{xx} \end{bmatrix}. 
\end{aligned}
\end{equation*}
We also apply a GARK method to this problem.  In this case, the right-hand side of ~\cref{eq:stiff_brusselator} is split into two operators:
\begin{equation*} \label{stiff_brus_gark}
  \begin{aligned}
    \Fl{1}(u,v,w) = \begin{bmatrix} \rho_u u_x \\ \rho_v v_x \\ \rho_w w_x \end{bmatrix},
    \Fl{2}(u,v,w) = \begin{bmatrix} \alpha_u u_{xx} + a - (w + 1)u + u^2v\\ \alpha_v v_{xx} + wu - u^2v\\ \alpha_w w_{xx} + (b-w)/\epsilon - wu \end{bmatrix}.
  \end{aligned}
  \end{equation*}
  The following multi-rate infinitesimal methods were tested: a
  third-order MRI-IMEX3 \cite{Chinomona2021}, a second-order MRI-IRK2
  and a third-order MRI-ESDIRK3a \cite{Sandu2019}.  The generalized
  additive Runge--Kutta method tested is a second-order
  stability-decoupled IMEX-GARK method, IMEX-GARK2, from example 6 of
  \cite{sandu2015} with $\beta = 1/2$. The second-order multi-rate
  method MrGARK-EX2-IM2 from \cite{sarshar2019} is also
  tested. Coefficients of these methods can be found in
  \cref{sec:os_methods_coeff}. We note that the MRI-IMEX3 method is
  applied to the $3$-split problem $\Fl{E} + \Fl{I} + \Fl{F}$. The
  MRI-IRK2, MRI-ESDIRK3a and MrGARK-EX2-IM2 methods are applied to the
  $2$-split problem $\Fl{S} + \Fl{F}$. The IMEX-GARK2 method is
  applied to the $2$-split problem $\Fl{1} + \Fl{2}$.  For the
  multi-rate infinitesimal methods, the adaptive solver used to find
  the exact solution of the fast operator is {\tt RK45} from \scipy\
  with \rtol=\num{1e-10} and \atol=\num{1e-12}.  The error is measured
  using the $\ell_2$-norm of the difference between the numerical and
  reference solutions over all the spatial grid points at $\tf = 3$.
  \Cref{fig:stiff_brusselator_convergence} shows that all five methods
  exhibit the expected order of convergence. Similar to observations
  in~\cite{Chinomona2021}, fourth-order methods experience order
  reduction for $\epsilon = 10^{-3}$ and have been omitted.

\section{Conclusions}
\label{sec:conclusions}

The \pythOS\ library has been designed to facilitate the systematic
numerical solution differential equations via operator splitting with
special emphasis on fractional-step methods.  A variety of classes of
operator-splitting methods have been implemented through similar
interfaces to allow for easy implementation and comparison between
methods, both within a class and between classes of methods.  The
solvers have been designed to be flexible, allowing various options
for sub-integrations where practical and any number of operators.  The
solvers have also been designed to solve ODEs or PDEs with support for
finite element discretization and complex inputs.  The examples
demonstrate the use of various non-standard options for
operator-splitting methods such as high-order methods, methods with
complex coefficients, and the use of arbitrary sub-integrators with
splittings of size $\Nop > 2$. Further examples that combine the easy
implementation of the hitherto under-explored GARK, multi-rate, and
multi-rate infinitesimal methods illustrate the utility of \pythOS\ as
a wide-ranging tool for exploration of operator-splitting methods.


\bibliographystyle{ACM-Reference-Format}
\bibliography{references}

\appendix

\section{Coefficients of operator-splitting methods}
\label{sec:os_methods_coeff}

\subsection{Fractional-step methods}

\begin{itemize}
\item The Lie--Trotter (Godunov) method is a first-order method that
  can be applied to an $\Nop$-additive system. Its coefficients
  are given in \cref{tab:godunov}.

\begin{table}[h!]
	\caption{Lie--Trotter coefficients}
	\label{tab:godunov}
\begin{tabular}{|c|c|}
	\hline 
	$k$ & $ \{\aaalpha{k}{\ell}\}_{\ell = 1,2,\dots, \Nop} $  \\
	\hline
	1 & 1  \\
	\hline
\end{tabular}
\end{table}

\item The Strang (Strang--Marchuk) method is a second-order method
  that can be applied to an $\Nop$-additive system. It can be
  seen as the composition of the Godunov method with its adjoint over
  half steps. Its coefficients are given in \cref{tab:strang}.

\begin{table}[h!]
\caption{Strang coefficients}
\label{tab:strang}
\begin{tabular}{|c|c|c|c|c|c|}
	\hline 
	$k$ & $ \aaalpha{k}{1} $  & $ \aaalpha{k}{2} $   & $\cdots$  & $ \aaalpha{k}{\Nop-1} $  & $ \aaalpha{k}{\Nop} $  \\
	\hline
	1 & $1/2$ &  $1/2$ &  $\cdots$ & $1/2$ & $1$ \\ 
	\hline 
	2 &  $0$ &  $\cdots$ &	$0$ 	&	$1/2$  & $0$ \\ 
	\hline
	$\vdots$ &  $\vdots$ &  $\vdots$ &  $\iddots$  &  $\vdots$  &  $\vdots$   \\ 
	\hline 
	$\Nop-1$ & $0$ & $1/2$ & $0$	&$\cdots$ & $0$ \\
	\hline 
	$\Nop$ & $1/2$ & $0$ & $\cdots$	&$\cdots$ & $0$ \\
	\hline 
\end{tabular}
\end{table}

\item PP3\_4A-3 \cite{auzinger2016practical} is a third-order
  real-valued method for $3$-additive system. Its coefficients
  are given in \cref{tab:pp34a3}.

\begin{table}[h!]
\caption{PP3\_4A-3 coefficients}
\label{tab:pp34a3}
\begin{tabular}{|c|r|r|r|}
	\hline 
	$k$ & \multicolumn{1}{|c|}{$ \aaalpha{k}{1} $}  & \multicolumn{1}{|c|}{$ \aaalpha{k}{2} $}  & \multicolumn{1}{|c|}{$ \aaalpha{k}{3}$} \\
	\hline
	1 & $0.461601939364879971$ & $-0.266589223588183997$ & $-0.360420727960349671$ \\ \hline
	2 & $-0.067871053050780081$ & $0.092457673314333835$ &  $0.579154058410941403$ \\ \hline
	3 & $-0.095886885226072025$ & $0.674131550273850162$ &  $0.483422668461380403$ \\ \hline
	4 & $0.483422668461380403$ &  $ 0.674131550273850162$ & $-0.095886885226072025 $ \\ \hline
	5 & $0.579154058410941403$ & $0.092457673314333835$ & $-0.067871053050780081 $ \\ \hline
	6 & $-0.360420727960349671$ & $-0.266589223588183997$ &  $0.461601939364879971$ \\ \hline
\end{tabular}
\end{table}

\item The Yoshida method is a fourth-order method that can be applied
  to an $\Nop$-additive split system. It is a three-fold composition
  of the Strang method over fractions of $\Delta t$. The coefficients
  of the Yoshida method are given in \cref{tab:yoshida}, where
  $\theta = 1/(2-\sqrt[3]{2})$.

\begin{table}[h!]
  \caption{Yoshida coefficients}
\label{tab:yoshida}
\begin{tabular}{|c|c|c|c|}
 	\hline 
	$k$ & $ \aaalpha{k}{1} $  & $ \aaalpha{k}{2} $   & $ \aaalpha{k}{3}$ \\
	\hline
	1 & $0$ & $0$ & ${\theta/2}$ \\  \hline
	2 & $0$ & ${\theta/2}$ & $0$ \\ \hline
	3 & $\theta$ & ${\theta/2}$ & ${(1-\theta)/2}$  \\ \hline
	4 & $0$ & ${(1-2\theta)/2}$ & $0$ \\ \hline
	5 & $1-2\theta$ & $(1-2\theta)/2$ & ${(1-\theta)/2}$ \\ \hline
	6 & $0$ & ${\theta/2}$ & $0$ \\ \hline
	7 & $\theta$ & $\theta/2$ & ${\theta/2}$ \\ \hline
\end{tabular}
\end{table}

\item The \CLTTwo\ method is a complex-valued, second-order method
  that can be generalized to $\Nop$ operators \cite{spiteri_wei_clt}. Its
  coefficients are given \cref{tab:clt2}.

\begin{table}[htbp]
\caption{Second-order complex Lie--Trotter}
\label{tab:clt2}
\begin{tabular}{|c|c|}
	\hline 
	k & $ \{\aaalpha{k}{\ell}\}_{\ell = 1,2,\dots, \Nop} $  \\
	\hline
	1 & $1/2+i/2$  \\ \hline 
	2 & $1/2-i/2$  \\
	\hline
\end{tabular}
\end{table}

\item The \CLTThree\ method is complex-valued, third-order method that
  can be generalized to $\Nop$ operators \cite{spiteri_wei_clt}.  The
  coefficients are given \cref{tab:clt3}.

\begin{table}[h!]
\caption{Third-order complex Lie--Trotter }
\label{tab:clt3}
\begin{tabular}{|c|c|}
	\hline 
	$k$ & $\displaystyle \{\alpha_k^{[\ell]}\}_{\ell = 1, 2,\dots, \Nop}$ \\
	\hline 
	1 & $\displaystyle 1/4- 1/(4\sqrt{3}) + \left(1/4 + 1/(4\sqrt{3}) \right)i$  \\
	\hline 
	2 &  $\displaystyle 1/4 + 1/(4\sqrt{3}) + \left(-1/4 + 1/(4\sqrt{3})\right)i$  \\
	\hline 
	3 & $\displaystyle 1/4 + 1/(4\sqrt{3}) + \left(1/4  - 1/(4\sqrt{3})\right)i$ \\
	\hline 
	4 &  $\displaystyle  1/4 - 1/(4\sqrt{3}) + \left(-1/4 - 1/(4\sqrt{3})\right)i$ \\
	\hline 
\end{tabular}
\end{table}

\end{itemize}

\subsection{GARK methods}

\begin{itemize}
\item IMEX-GARK2 is a second-order, stability-decoupled IMEX-GARK
  method with $\beta = 1/2$ \cite{sandu2015}. The Butcher tableau
  representation is given below.

$\begin{array}{c|c|c|c}
	  \cc{E,E}   &   \AAA{E,E}   &   \AAA{E,I}   &   \cc{E,I}   \\ \hline
	    \cc{I,E}   &   \AAA{I,E}   &   \AAA{I,I}   &   \cc{I,I}   \\ \hline 
	  &   \bb{E}    &    \bb{I}    & \\
\end{array} := 
\begin{array}{c|c|c|c}
	\begin{array}{c}
		0 \\ 
		1/2 \\ 
		1 \\
	\end{array} & 
\begin{array}{ccc}
	0 & 0 & 0 \\
	1/2 & 0 & 0  \\ 
	1-\beta & \beta & 0 \\
\end{array} & 
\begin{array}{cc}
	0 & 0 \\
	1/2 & 0 \\
	1/2 & 1/2 \\ 
\end{array} & 
\begin{array}{c}
	 0 \\ 
	 1/2 \\ 
	 1 \\
\end{array} \\ \hline 
\begin{array}{c}
 1/4 \\ 
 3/4 \\ 
\end{array}& 
\begin{array}{ccc}
1/4 & 0 & 0 \\
1/4 & 1/2 & 0 \\ 
\end{array}
& 
\begin{array}{cc}
1/4 & 0 \\ 
1/2 & 1/4 \\ 
\end{array}
& 
\begin{array}{c}
	1/4 \\ 
	3/4 \\ 
\end{array} \\ 
\hline 
& 
\begin{array}{ccc}
	1/4 & 1/2 & 1/4 \\
\end{array}
& 
\begin{array}{cc}
	1/2  & 1/2 \\
\end{array} & 
\end{array}$

\end{itemize}

\subsection{Multi-rate methods}

\begin{itemize}
\item MRI-IRK2 is a second-order MRI-GARK method based on the implicit
  trapezoidal rule \cite{Sandu2019}. The coefficients are given below.

 	$$\cc{S} = \begin{bmatrix}
		0 \\ 
		1 \\ 
		1 \\ 
		1 \\
      \end{bmatrix},\
	\GGamma{0} = \begin{bmatrix}
		0 & 0 & 0 & 0 \\ 
		1 & 0 & 0 & 0 \\ 
		-1/2 & 0 & 1/2 & 0 \\ 
		0 & 0 & 0 & 0 \\
	\end{bmatrix}. 
	$$

  \item MRI-ESDIRK3a is a third-order, stiffly accurate method with
    equidistant abscissae \cite{Sandu2019}. Let
    $\lambda = 0.435866521508458999416019$, the coefficients are given
    below.
	
	
	$$\cc{S} = \begin{bmatrix}
		0 \\ 1/3 \\ 1/3 \\ 2/3 \\ 2/3 \\ 1 \\ 1 \\ 1 \\ 
	\end{bmatrix},$$ 

	$\GGamma{0} = \begin{bmatrix}
		 0 & 0 & 0 & 0 & 0 & 0 & 0 & 0 \\
		 1/3 & 0 & 0 & 0 & 0 & 0 & 0 & 0 \\
		 - \lambda &  0 &   \lambda &  0 &  0 &  0 &  0 &  0 \\
		 (3-10  \lambda)/(24  \lambda - 6) &  0 &  (5-18  \lambda)/(6-24  \lambda) &  0 &  0 &  0 &  0 &  0 \\
		 (-24  \lambda^2 + 6  \lambda+1)/(6-24  \lambda) &  0 &  (-48  \lambda^2+12  \lambda+1)/(24  \lambda-6) &  0 &   \lambda &  0 &  0 &  0 \\
		 (3-16  \lambda)/(12-48  \lambda) &  0 & (48  \lambda^2 - 21  \lambda + 2)/(12  \lambda-3) & 0 & (3-16  \lambda)/4 &  0 &  0 &  0 \\
		 - \lambda &  0 &  0 &  0 &  0 &  0 &   \lambda &  0 \\
		 0 & 0 & 0 & 0 & 0 & 0 & 0 & 0 \\ 
	\end{bmatrix}$

  \item MRI-IMEX3 is a third-order MRI-IMEX method
    \cite{Chinomona2021}. The non-zero coefficients (accurate to 36
    decimal places) are
	
	$\cc{S}_2 =  \cc{S}_3 =  0.4358665215084589994160194511935568425 $ \\ 
	
	$\cc{S}_4 =  \cc{S}_5 = 	0.7179332607542294997080097255967784213$ \\ 
	
	$\cc{S}_6 =  \cc{S}_7 =  \cc{S}_8 = 1$ \\ 

	$\GGamma{0}_{2 \\1} = \GGamma{0}_{3 \\1} = \GGamma{0}_{3 \\3} = \GGamma{0}_{5 \\5} = \GGamma{0}_{6 \\1} = -\GGamma{0}_{7 \\1} = \GGamma{0}_{7 \\7} = 0.4358665215084589994160194511935568425$ \\ 
	
	$\GGamma{0}_{4 \\1} = - \GGamma{0}_{5 \\1} = - 0.4103336962288525014599513720161078937$ \\ 
	
	$\GGamma{0}_{4 \\3} = 0.6924004354746230017519416464193294724$ \\ 
	
	$\GGamma{0}_{5 \\3} =  -0.8462002177373115008759708232096647362$ \\ 
	
	$\GGamma{0}_{6 \\3} = 0.9264299099302395700444874096601015328$ \\ 
	
	$\GGamma{0}_{6 \\5} = -1.080229692192928069168516586450436797$ \\

	$\OOmega{0}_{2 \\1} = \OOmega{0}_{8 \\7} = 0.4358665215084589994160194511935568425$ \\
	
	$\OOmega{0}_{4 \\1} = -0.5688715801234400928465032925317932021$ \\
	
	$\OOmega{0}_{4 \\3} = 0.8509383193692105931384935669350147809$ \\
	
	$\OOmega{0}_{5 \\1} = -\OOmega{0}_{5 \\3} = 0.454283944643608855878770886900124654$ \\
	
	$\OOmega{0}_{6 \\1} = - 0.4271371821005074011706645050390732474$ \\ 
	
	$\OOmega{0}_{6 \\3} = 0.1562747733103380821014660497037023496$ \\
	
	$\OOmega{0}_{6 \\5} = 0.5529291480359398193611887297385924765$ \\
	
	$\OOmega{0}_{8 \\1} = 0.105858296071879638722377459477184953$ \\
	
	$\OOmega{0}_{8 \\3} = 0.655567501140070250975288954324730635$ \\
	
	$\OOmega{0}_{8 \\5} = -1.197292318720408889113685864995472431$.

  \item MrGARK EX2-IM2 is a second-order MrGARK method with an
    explicit fast part and an implicit slow
    part~\cite{sarshar2019}. The coefficients are given below.
	
	$\AAA{F,F} = \begin{bmatrix}
		0 & 0 \\
		2/3 & 0 \\
	\end{bmatrix}$, 
	$\AAA{S,S} = \begin{bmatrix}
	1-1/\sqrt{2} & 0 \\
	1/\sqrt{2} & 1-1/\sqrt{2} \\
	\end{bmatrix}$, 
	$\AAA{S,F,1} = \begin{bmatrix}
		M-M/\sqrt{2} & 0 \\
		1/4 & 3/4 \\
	\end{bmatrix}$, 
	
	$\AAA{F,S,\faststepind} = \begin{bmatrix}
		(\faststepind-1)/M & 0 \\
		(3\faststepind-1)/(3M) & 0 \\
	\end{bmatrix}$,\ $\faststepind = 1,2,\dots, M$, 
	
	$\AAA{S,F,\faststepind} = \begin{bmatrix}
		0 & 0 \\
		1/4 & 3/4 \\
	\end{bmatrix}$,\ $\faststepind = 2,3,\dots, M$,
	
	$\bb{F,F} = \begin{bmatrix} 
		1/4 & 3/4 \\ 
		\end{bmatrix}
		$,	
	$\bb{S,S} = \begin{bmatrix} 
		1/\sqrt{2} & 1-1/\sqrt{2}\\ 
	\end{bmatrix}.
	$
	
\end{itemize}





\end{document}